\documentclass[epjST]{svjour}%
\usepackage{amsfonts, amsmath, amssymb}
\usepackage{bmpsize}
\usepackage{bmpsize}
\usepackage{graphics}
\usepackage{amsmath}
\usepackage{amsfonts}
\usepackage{amssymb}
\usepackage{graphicx}%
\providecommand{\U}[1]{\protect\rule{.1in}{.1in}}
\begin{document}
\title{Symmetry-breaking for a restricted $\pmb{n}$-body problem in the Maxwell-ring configuration}
\author{
Renato Calleja \inst{1}
\fnmsep\thanks
{\email{calleja@mym.iimas.unam.mx}}
\and
Eusebius Doedel \inst{2}
\fnmsep\thanks
{\email{doedel@cs.concordia.ca}}
\and
Carlos Garc\'{\i}a-Azpeitia \inst{3}
\fnmsep\thanks
{\email{cgazpe@ciencias.unam.mx}}
}

\institute{Matem\'aticas y Mec\'anica, IIMAS, Universidad Nacional Aut\'onoma de
M\'exico, Admon. No. 20, Delegaci\'on Alvaro Obreg\'on, 01000 M\'exico D.F.
\and
Department of Computer Science, Concordia University, 1455 boulevard de
Maisonneuve O., Montr\'eal, Qu\'ebec H3G~1M8, Canada \and
Departamento de Matem\'{a}ticas, Facultad de Ciencias, Universidad Nacional
Aut\'{o}noma de M\'{e}xico, 04510 M\'{e}xico DF, Mexico}




\abstract{
We investigate the motion of a massless body interacting with the Maxwell
relative equilibrium, which consists of $n$ bodies of equal mass at the
vertices of a regular polygon that rotates around a central mass.
The massless body has three equilibrium $\mathbb{Z}_{n}$-orbits from which
families of Lyapunov orbits emerge.
Numerical continuation of these families using a boundary value formulation
is used to construct the bifurcation diagram for the case $n=7$, also
including some secondary and tertiary bifurcating families.
We observe symmetry-breaking bifurcations in this system, as well as certain
period-doubling bifurcations.
}
\maketitle

\section*{Introduction}

In his 1859 essay \cite{Ma} Maxwell proposed a model to study the rings of
Saturn. His model consists of $n$ bodies of equal mass at the vertices of a
regular polygon that rotates around a massive body at the center. Maxwell used
Fourier analysis and dispersion relations in the determination of the
stability of the ring. The Maxwell equilibrium has been studied in several
papers since then. In particular, Moeckel proved in \cite{Moekel92} that the
equilibrium is stable if $n\geq7$ and the body at the center massive enough.
See also \cite{GaIz13,VaKo07,Ro00} and references therein.

In this paper we consider the motion of a satellite under the gravitational
effect of the Maxwell equilibrium. Several papers have been devoted to study
the stability and bifurcation of periodic solutions for the restricted
$N$-body problem in the Maxwell configuration. For example, a study of the
existence and linear stability of equilibrium positions can be found in
\cite{BaEl04}, an analysis of the bifurcation of planar and vertical families
of periodic solutions in \cite{GaIz10}, and a numerical exploration in
\cite{Ka08}.

Given that a change of stability occurs at $\mu_{1} \approx584$ when $n=7$, we
present a numerical exploration of the motion of the satellite for this stable
system, taking $\mu=10^{3}$. We follow the planar and vertical Lyapunov
families that were proved to exist in \cite{GaIz10}, and we compute new,
secondary families that bifurcate from the Lyapunov families. These numerical
results allow us to construct a bifurcation diagram in Figure~\ref{fig:2} that
shows primary, secondary, and tertiary bifurcating families. The results also
allow us to present the isotropy lattice in Figure~\ref{fig:3} which shows the
symmetries of these families.

The numerical computations in this paper are done using continuation methods
and boundary value techniques for determining the periodic orbits that emanate
from the equilibrium orbits. Python scripts that make the AUTO software
perform the calculations reported in this paper will be made freely available.
Similar techniques have been applied to the restricted $3$-body problem; see
for example \cite{DoVa}, where a detailed bifurcation diagram with various
families of periodic orbits can be found.

This paper is organized as follows. In Section~\ref{sec:restrictedMaxwell} we
recall some key results from the literature concerning the equilibria and the
Lyapunov orbits of the problem. In Section~\ref{sec:symmetryBreaking} we
present a bifurcation diagram and an isotropy lattice for the restricted
$N$-body problem in the Maxwell configuration, with $n=7$ and $\mu=10^{3}$. In
Section~\ref{sec:first} we describe the isotropy groups of the Lyapunov
families. In Section~\ref{sec:second} we address secondary bifurcations, and
in Section~\ref{sec:third} describe some tertiary families. We also present
evidence of the existence of invariant tori foliated by periodic orbits.
Finally, in Section~\ref{sec:four}, we consider the breaking of symmetries of
planar interplanetary periodic orbits.

\section{The restricted $N$-body problem}

\label{sec:restrictedMaxwell}

The Maxwell relative equilibrium consists of a body of mass $\mu$ at
$\pmb{a}_{0}=\mathbf{0\in\mathbb{R}^{3}}$, and $n$ bodies of mass $1$ located
at $\pmb{a}_{j}=(e^{ij\zeta},0)\in\mathbb{R}^{3}$, $j=1,\cdots,n$, where
$\zeta=2\pi/n$. These positions and masses correspond to a relative
equilibrium solution of Newton's equations, $\pmb{q}_{j}(t)=e^{\mathcal{J}
t}\pmb{a}_{j}$, when the masses are renormalized by $\mu+s$
\cite{BaEl04,GaIz10}, where $\mathcal{J}=diag(J,0)$,
\[
J=\left(
\begin{array}
[c]{cc}%
0 & -1\\
1 & 0
\end{array}
\right)  ,\qquad s=\frac{1}{4}\sum_{j=1}^{n-1}\frac{1}{\sin(j\pi/n)}\text{.}
\]

The equation of a satellite in rotating coordinates,
$\pmb{q}(t)=e^{\mathcal{J}t}\pmb{u}(t)$, is
\begin{equation}
\ddot{\pmb{u}}+2\mathcal{J}\dot{\pmb{u}}=\nabla V\text{,}\label{satellite}%
\end{equation}
where $\pmb{u}=\left(  x,y,z\right)  \in\mathbb{R}^{3}$ and
\begin{equation}
V(\pmb{u})=\frac{1}{2}\left\Vert (x,y)\right\Vert ^{2}+\frac{\mu}{s+\mu}%
\frac{1}{\left\Vert \pmb{u}\right\Vert }+\sum_{j=1}^{n}\frac{1}{s+\mu}\frac
{1}{\left\Vert \pmb{u}-(e^{ij\zeta},0)\right\Vert }\text{.}\label{potential}%
\end{equation}
The first term of the potential $V(\pmb{u})$ in Equation~\eqref{potential}
corresponds to the centrifugal force, the second term is the interaction with
the mass $\mu$, and the third term models the interaction with the $n$
primaries of mass $1$; see \cite{BaEl04} and \cite{GaIz10}.
\begin{figure}[t]
\begin{center}
\resizebox{13cm}{!} {\includegraphics{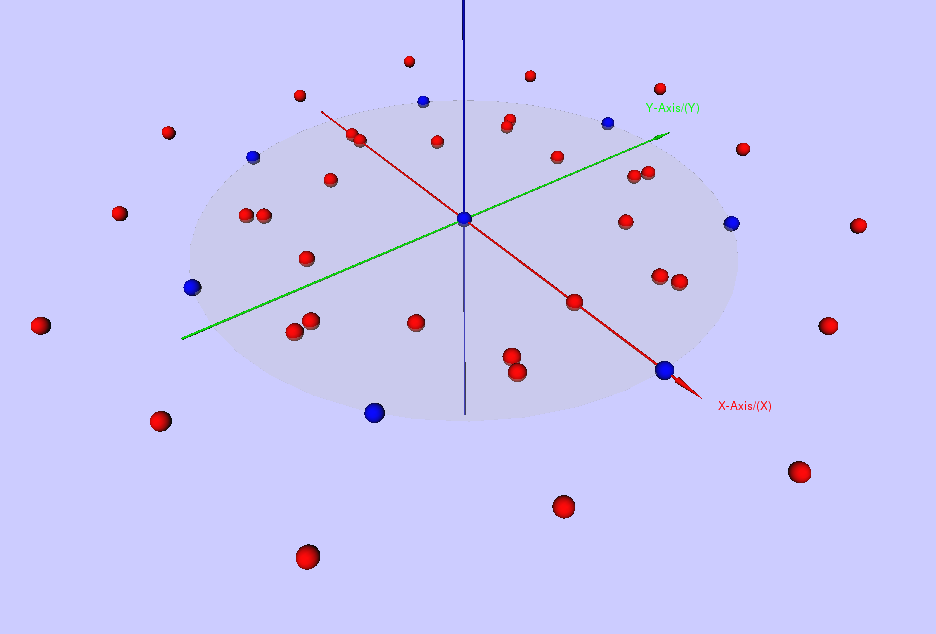}  \hskip1.0cm
\includegraphics{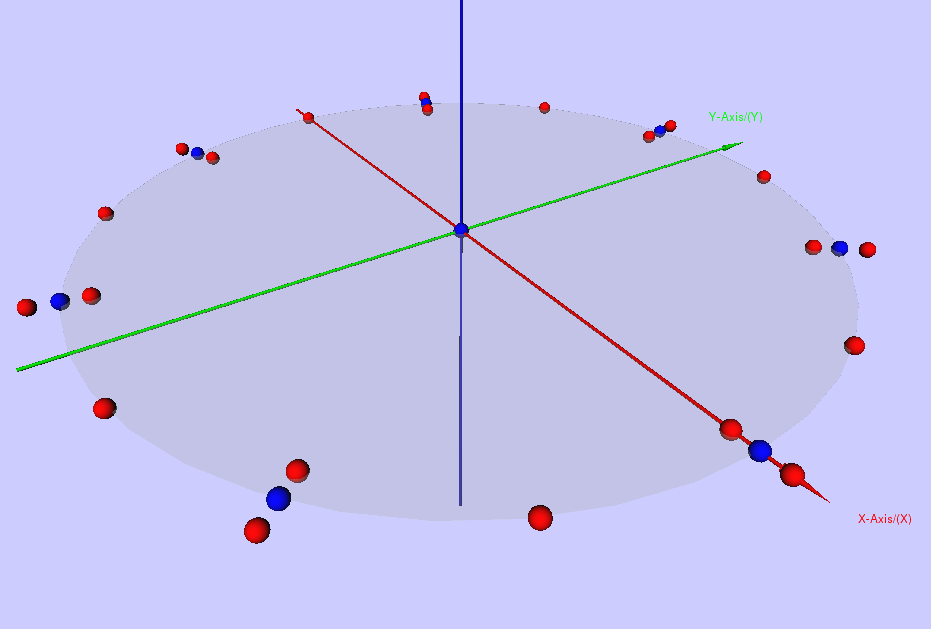} }
\end{center}
\caption{The Maxwell configuration with the masses colored blue and the
equilibrium orbits colored red. Left: the $35$ equilibria when $\mu=3$. Right:
the $21$ equilibria when $\mu=1000$. }%
\label{fig:1}%
\end{figure}
The equilibria of equation (\ref{satellite}) are critical points of $V$,
defined in Equation~\eqref{potential}. Moreover, due to the particular form of
the Maxwell configuration, the potential $V$ is $\mathbb{Z}_{n}$-invariant. In
this respect the existence of three $\mathbb{Z}_{n}$-orbits of equilibria for
any value of $\mu$ is proved in \cite{BaEl04,GaIz10}. For small $\mu$, two
additional $\mathbb{Z}_{n}$-orbits of equilibria appear close to the origin
\cite{GaIz10}; see Figure~\ref{fig:1}.

In \cite{GaIz10}, taking advantage of the symmetries of the equations, the
authors analyze the bifurcation of periodic solutions from the $\mathbb{Z}
_{n}$-orbits of equilibria. Here we state these results for the equilibria
$\mathcal{L}_{1}$, $\mathcal{L}_{2}$ and $\mathcal{L}_{3}$ of the three
$\mathbb{Z}_{n}$-orbits:

\newpage

\begin{theorem}
\textrm{[Ize \& Garc\'{\i}a-Azpeitia \cite{GaIz10}.]} The libration point
$\mathcal{L}_{1}$, has one global bifurcation family of planar periodic
solutions that will be denoted by $L_{1}$, and one global family of vertical
solutions, denoted by $V_{1}$. Similarly, the libration point $\mathcal{L}
_{2}$ has one global family of planar periodic solutions, denoted $L_{2}$, and
one family of vertical solutions, $V_{2}$. For $\mu>\mu_{1}$, the equilibrium
$\mathcal{L}_{3}$ has two global bifurcations of planar solutions, one of
which has longer period, denoted $L_{3l}$, and another planar family of
shorter period, $L_{3s}$. There is also a bifurcating family of vertical
periodic solutions that will be denoted by $V_{3}$. Moreover, as a consequence
of the symmetries, the shape of all vertical solutions close to the
equilibrium resembles a spatial figure eight.
\end{theorem}

The \textit{global} property guarantees that the family is a continuum that
either goes to infinity in Sobolev norm or period, ends in a collision, or
ends at a bifurcation point. Indeed, each one of these possibilities appears
in the numerical continuation of the families, as illustrated in the
bifurcation diagram in Figure~\ref{fig:2}.

\section{Breaking of symmetries}

\label{sec:symmetryBreaking}

In this section we discuss the breaking of symmetries of
Equation~\eqref{satellite} for the case $n=$ $7$ and $\mu=10^{3}$, as observed
in the numerically computed Lyapunov families that emerge from the libration
points and the secondary families that bifurcate from them.

The $2\pi/\nu$-periodic solutions of Equation~(\ref{satellite}) are zeros of
the map
\[
\pmb{f}(\pmb{u};\nu) = \nu^{2}\ddot{\pmb{u}}+\mathcal{J}\nu\dot{\pmb{u}}
-\nabla V(\pmb{u})\text{,}
\]
defined in a set of $2\pi$-periodic collisionless functions $\pmb{u}$; see
\cite{GaIz10}. Since the potential is $\mathbb{Z}_{7}$-invariant and the
equations are autonomous, the map $\pmb{f}$ is equivariant under the action of
$(\zeta,\varphi)\in\mathbb{Z}_{7}\times S^{1}$ given by
\[
\rho(\zeta,\varphi)\pmb{u}(t)=e^{\mathcal{J}\zeta}\pmb{u}(t+\varphi)\text{,}
\]
where $\zeta=2\pi/7$. In addition the equations are symmetric with respect to
reflection of $y$ about the $xz$ plane, while reversing time, and with respect
to reflection of $z$ about the $xy$ plane. In this regard we define the
reflections $\kappa_{y}$ and $\kappa_{z}$ by
\[
\rho(\kappa_{y})\pmb{u}(t)=R_{y}\pmb{u}(-t)\text{~}~\mathrm{and}~~ \rho
(\kappa_{z})\pmb{u}(t)=R_{z} \pmb{u}(t)\text{,}
\]
where $R_{y}=diag(1,-1,1)$ and $R_{z}=diag(1,1,-1)$. Therefore the map
$\pmb{f}$ is equivariant under the full symmetry group
\begin{equation}
G=\left(  \mathbb{Z}_{7}\times S^{1}\cup\kappa_{y}(\mathbb{Z}_{7}\times
S^{1})\right)  \times\mathbb{Z}_{2}(\kappa_{z}). \label{fullsymgroup}%
\end{equation}

We will use the property that the group orbit of a function $\pmb{u}$,
\[
G(\pmb{u})=\{\rho(\gamma)\pmb{u}:\gamma\in G\}\text{,}%
\]
is isomorphic to $G/G_{\pmb{u}}$, where $G_{\pmb{u}}$ is the isotropy group
defined as
\[
G_{\pmb{u}}=\{\pmb{u}:\rho(\gamma)\pmb{u}=\pmb{u},\forall\gamma\in G\}\text{.}%
\]
For example, the libration equilibria $\mathcal{L}_{j}$, for $j=1,2,3$, have
group orbits $G(\mathcal{L}_{j})\simeq\mathbb{Z}_{7}$ and isotropy group
\[
G_{\mathcal{L}_{j}}=(S^{1}\cup\kappa_{y}S^{1})\times\mathbb{Z}_{2}(\kappa
_{z})\text{.}%
\]
We therefore present the breaking of symmetries only for the libration points
$\mathcal{L}_{j}$, $j=1,2,3$. \vskip0.25cm
\begin{figure}[t]
\begin{center}
\resizebox{13.2cm}{!}{\includegraphics{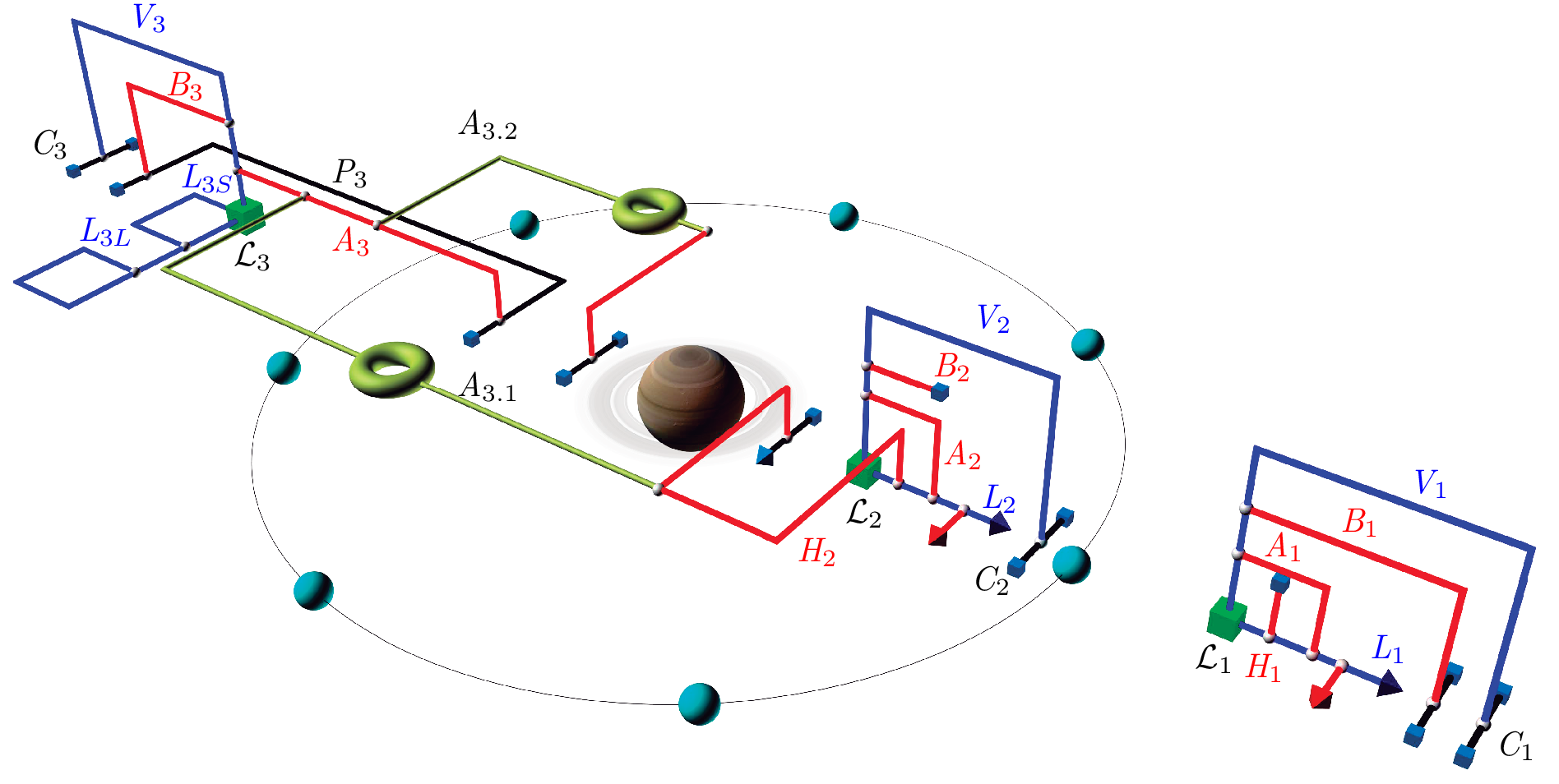} }
\end{center}
\caption{ Bifurcation diagram for $\mu=10^{3}$: The green cubes denote the
equilibrium positions $\mathcal{L}_{1}$, $\mathcal{L}_{2}$, and $\mathcal{L}%
_{3}$, the white spheres represent bifurcation points, the blue cubes
represent families that apparently end in collisions, and tetrahedra indicate
that the family goes to infinity in period or in Sobolev norm.}%
\label{fig:2}%
\end{figure}

The bifurcation diagram for the case $\mu=10^{3}$ is given in
Figure~\ref{fig:2}. The blue lines represent the vertical Lyapunov families
$V_{j}$ and the planar Lyapunov families $L_{j}$, $j=1,2,3$. Planar families
are positioned in the plane of the bodies, as are the black lines that
represent planar interplanetary orbits. The red lines are the result of a
secondary symmetry-breaking, with two solutions per equilibrium, as for the
Halo orbits $H_{j}$ and the Axial orbits $A_{j}$, $j=1,2,3$. The green lines
correspond to a tertiary symmetry-breaking, and as such they have a trivial
isotropy group and four symmetry-related branches per libration point.

\begin{remark}
: Some of the families in the bifurcation diagram that end at a tetrahedron,
in fact terminate as a heteroclinic orbit. For the restricted three-body
problem similar heteroclinic orbits are given in \cite{CaDo12,LoMa00}, and the
existence of heteroclinic connections is proved in \cite{LiMa}. In future work
we will present many other heteroclinic connections that we have located by
continuation of orbits in stable/unstable manifolds. Furthermore, in future
work we will present evidence of families of planar orbits that interconnect
the planar families $L_{3s}$ and $L_{3l}$, as in \cite{Henrard}. All results
are accompanied by scripts that allow their reproduction.
\end{remark}

\begin{figure}[th]
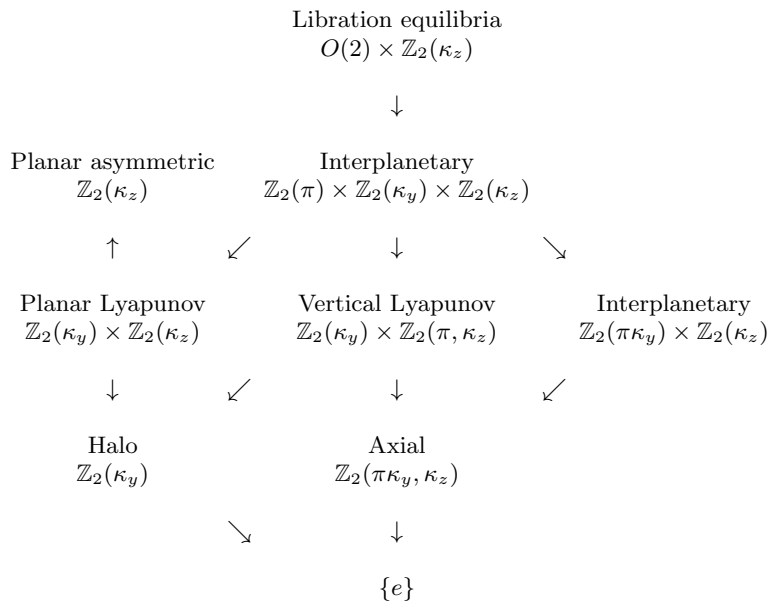

\centering
\par%
\[%
\begin{array}
[c]{ccccc}
&  &
\begin{array}
[c]{c}%
\text{Libration equilibria}\\
O(2)\times\mathbb{Z}_{2}(\kappa_{z})
\end{array}
&  & \\
&  &  &  & \\
&  & \downarrow &  & \\
&  &  &  & \\%
\begin{array}
[c]{c}%
\text{Planar asymmetric}\\
\mathbb{Z}_{2}(\kappa_{z})
\end{array}
&  &
\begin{array}
[c]{c}%
\text{Interplanetary}\\
\mathbb{Z}_{2}(\pi)\times\mathbb{Z}_{2}(\kappa_{y})\times\mathbb{Z}_{2}
(\kappa_{z})
\end{array}
&  & \\
&  &  &  & \\
\uparrow & \swarrow & \downarrow & \searrow & \\
&  &  &  & \\%
\begin{array}
[c]{c}%
\text{Planar Lyapunov}\\
\mathbb{Z}_{2}(\kappa_{y})\times\mathbb{Z}_{2}(\kappa_{z})
\end{array}
&  &
\begin{array}
[c]{c}%
\text{Vertical Lyapunov}\\
\mathbb{Z}_{2}(\kappa_{y})\times\mathbb{Z}_{2}(\pi,\kappa_{z})
\end{array}
&  &
\begin{array}
[c]{c}%
\text{Interplanetary}\\
\mathbb{Z}_{2}(\pi\kappa_{y})\times\mathbb{Z}_{2}(\kappa_{z})
\end{array}
\\
&  &  &  & \\
\downarrow & \swarrow & \downarrow & \swarrow & \\
&  &  &  & \\%
\begin{array}
[c]{c}%
\text{Halo}\\
\mathbb{Z}_{2}(\kappa_{y})
\end{array}
&  &
\begin{array}
[c]{c}%
\text{Axial}\\
\mathbb{Z}_{2}(\pi\kappa_{y},\kappa_{z})
\end{array}
&  & \\
&  &  &  & \\
& \searrow & \downarrow &  & \\
&  &  &  & \\
&  & \{e\} &  &
\end{array}
\]
\caption{Isotropy lattice}%
\label{fig:3}%
\end{figure}

The breaking of symmetries in the bifurcation diagram gives rise to the
lattice of isotropy groups of bifurcating orbits (isotropy lattice) in
Figure~\ref{fig:3}.

\section{Lyapunov families}

\label{sec:first}

The first bifurcation occurs when the $S^{1}$-symmetry is broken, giving rise
to the Lyapunov families (the blue lines in Figure~\ref{fig:2}) from the
equilibria $\mathcal{L}_{j}$, $j=1,2,3$. The number of symmetry-related
Lyapunov branches is equal to the order of $G/G_{\pmb{u}}$, which is $7$,
since each isotropy group is isomorphic to $\mathbb{Z}_{2}\times\mathbb{Z}%
_{2}$.
\begin{figure}[t]
\begin{center}
\resizebox{13cm}{!} {\includegraphics{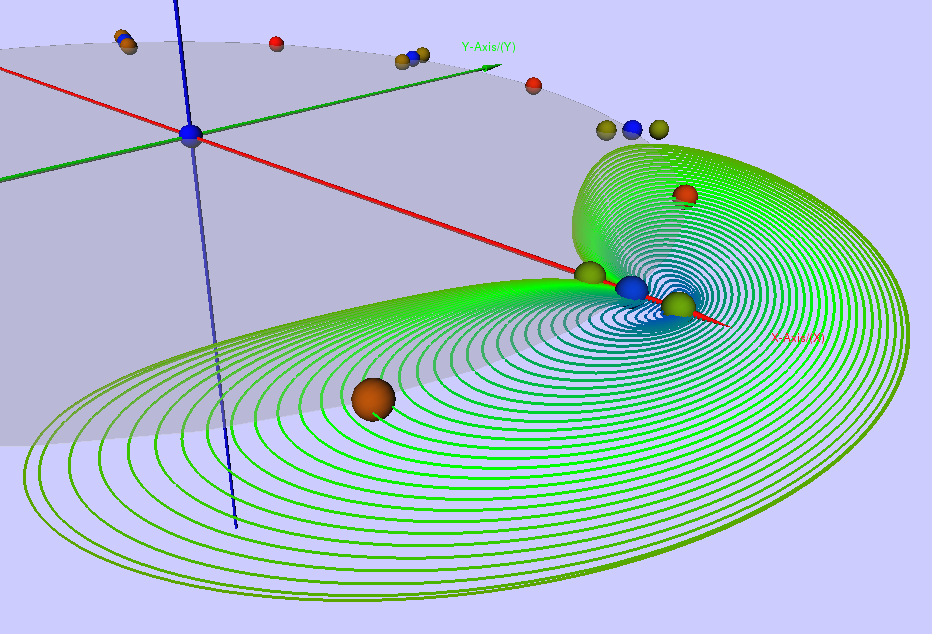}  \hskip1.9cm
\includegraphics{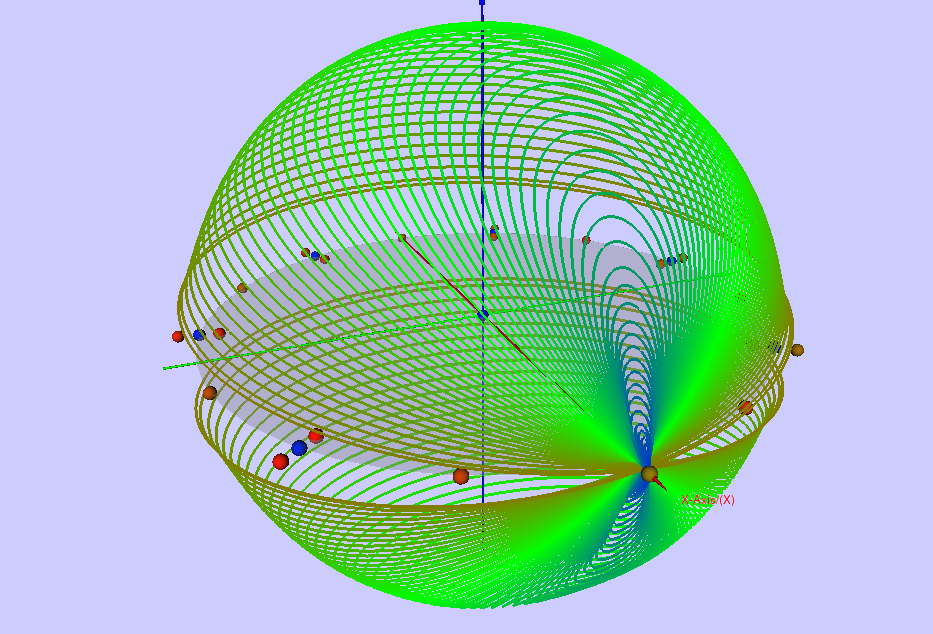} } \vskip0.5cm
\resizebox{13cm}{!} {\includegraphics{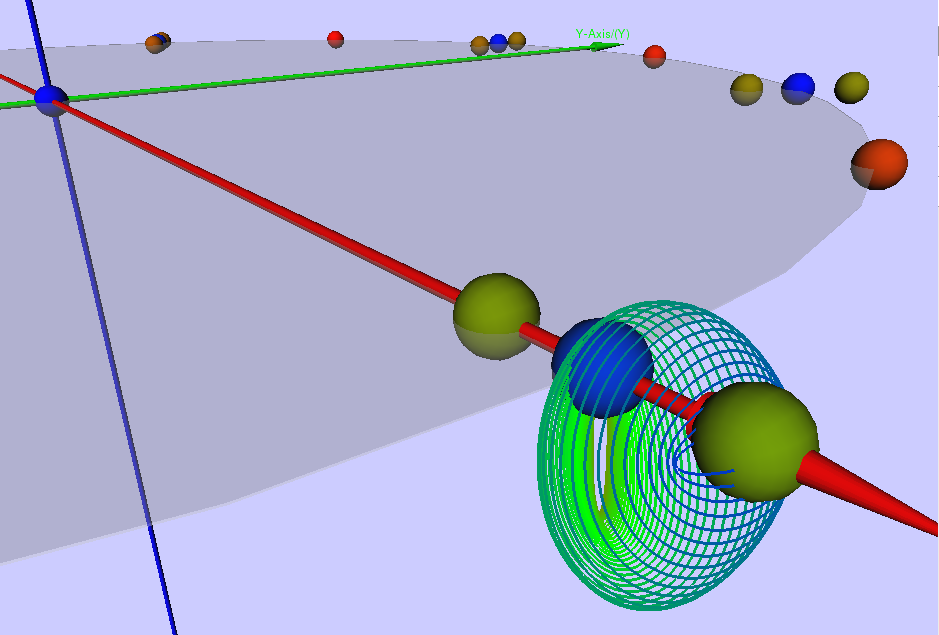}  \hskip1.9cm
\includegraphics{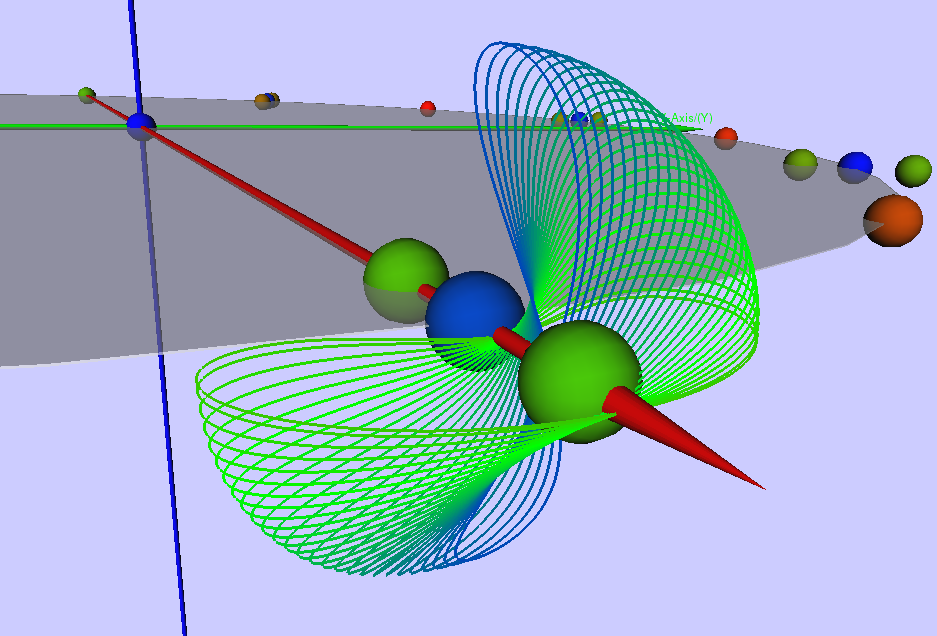} } \vskip0.5cm
\resizebox{13cm}{!} {\includegraphics{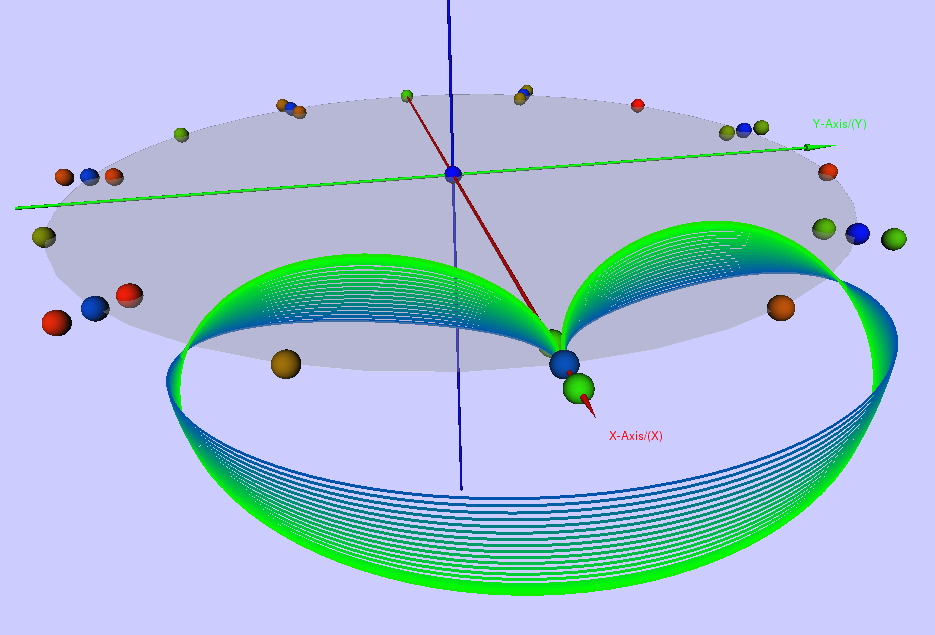}  \hskip1.9cm
\includegraphics{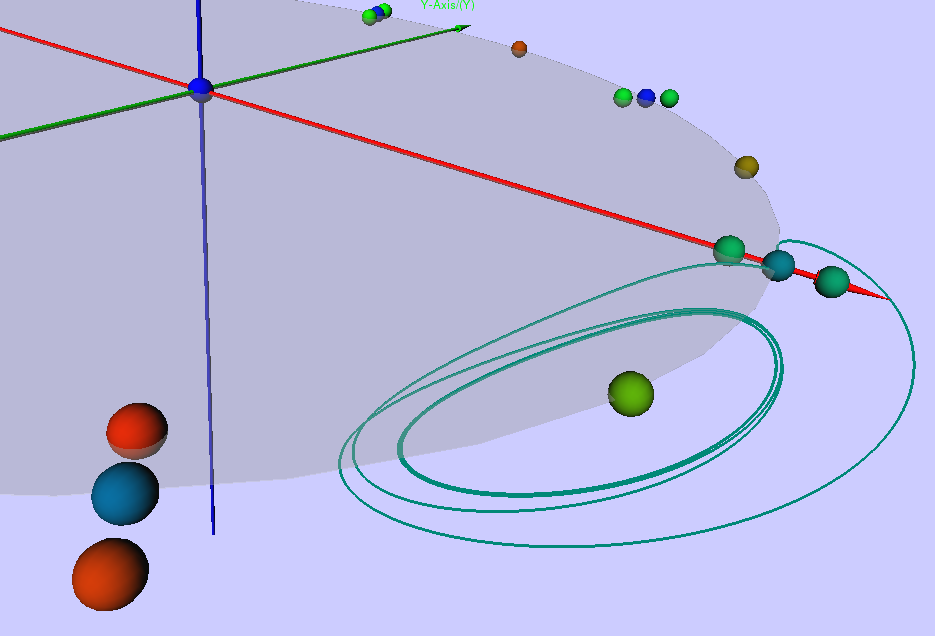} }
\end{center}
\caption{ Top Left: The planar Lyapunov family $L_{1}$, which ends in a
collision orbit. Top Right: The vertical family $V_{1}$ until its second
bifurcation orbit. Center Left: The \textquotedblleft
Southern\textquotedblright\ branch of the Halo family $H_{1}$, which ends in a
collision orbit. Center Right: One branch of the Axial family $A_{1}$, which
forms a \textquotedblleft bridge\textquotedblright\ between $L_{1}$ and
$V_{1}$. Bottom Left: One branch of the third family that bifurcates from
$L_{1}$, which ends in a collision orbit. Bottom Right: One branch of the
fourth family that bifurcates from $L_{1}$, which approaches an orbit that is
homoclinic to a periodic orbit. }%
\label{fig:4}%
\end{figure}

Hereafter $H<G$ denotes that $H$ is a subgroup of $G$. The isotropy group of
the planar Lyapunov orbits that emerge from the libration equilibria is
\begin{equation}
G_{L_{j}}=\mathbb{Z}_{2}(\kappa_{y})\times\mathbb{Z}_{2}(\kappa_{z}%
)<G_{\mathcal{L}_{j}}\text{.}\label{GP}%
\end{equation}
These periodic solutions have the property that $x(t)$ is even, $y(t)$ is odd,
and $z(t)=0$. In particular, we observe that the planar orbits are invariant
under the transformation that takes $y$ to $-y$.

For the vertical Lyapunov orbits the isotropy group is
\begin{equation}
G_{V_{j}}=\mathbb{Z}_{2}(\kappa_{y})\times\mathbb{Z}_{2}(\pi,\kappa
_{z})<G_{\mathcal{L}_{j}}\text{.} \label{GV}%
\end{equation}
Here a solution is fixed by $G_{V_{j}}$ if it satisfies
\[
\pmb{u}(t)=\rho(\kappa,\pi)\pmb{u}(t)=R_{z}\pmb{u}(t+\pi)\text{,}
\]
which is equivalent to assuming that $x(t)$ is a $\pi$-periodic even function,
$y(t)$ is a $\pi$-periodic odd function, and $z(t)=-z(t+\pi)$. Therefore these
solutions follow the planar $\pi$-periodic curve $(x,y)$ twice; one time with
the spatial coordinate $z$ and a second time with $-z$. This fact was proved
in \cite{GaIz10}.
\begin{figure}[t]
\begin{center}
\resizebox{13cm}{!} {\includegraphics{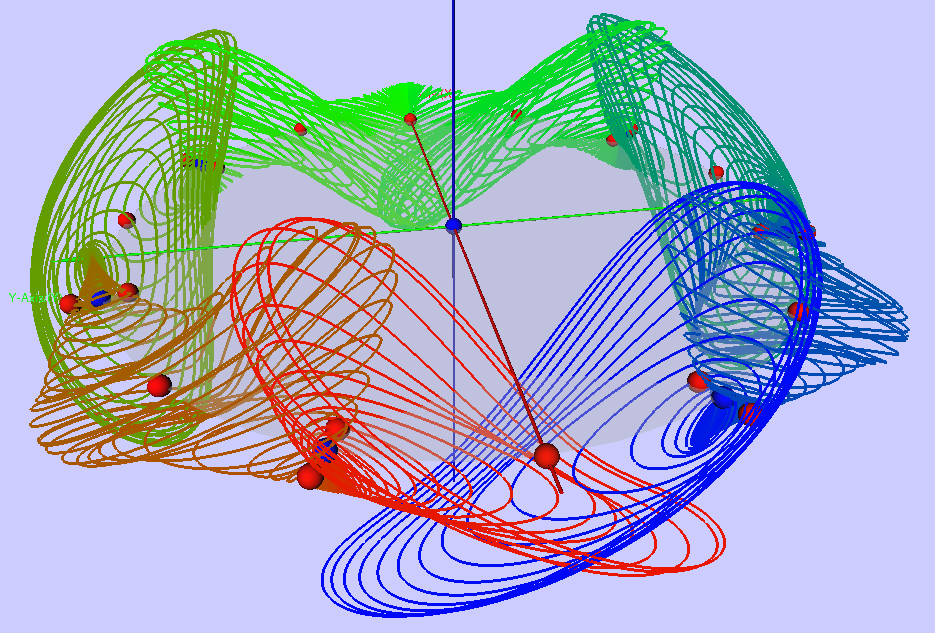}  \hskip1.9cm
\includegraphics{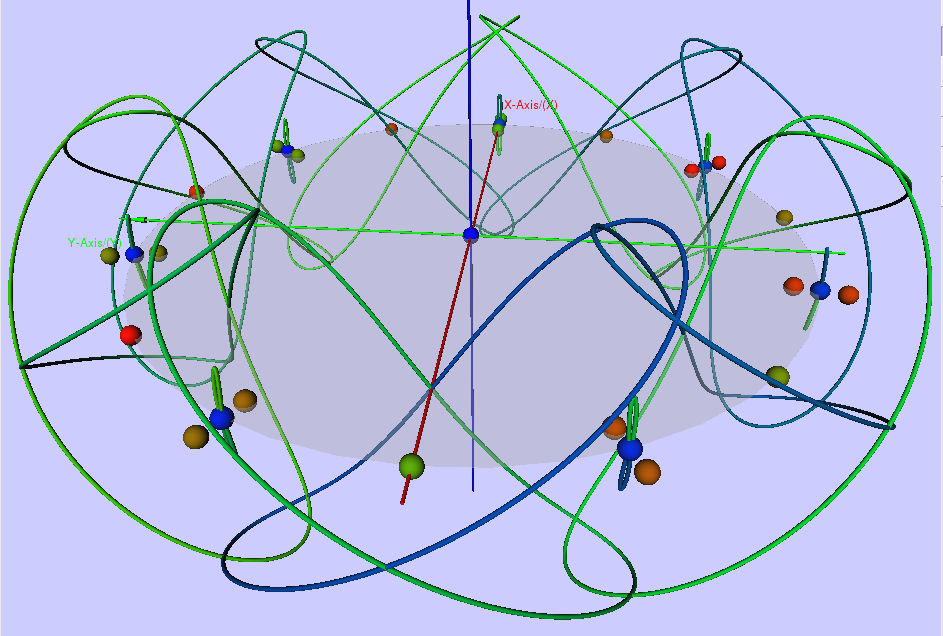} }
\par
\resizebox{13cm}{!} {\includegraphics{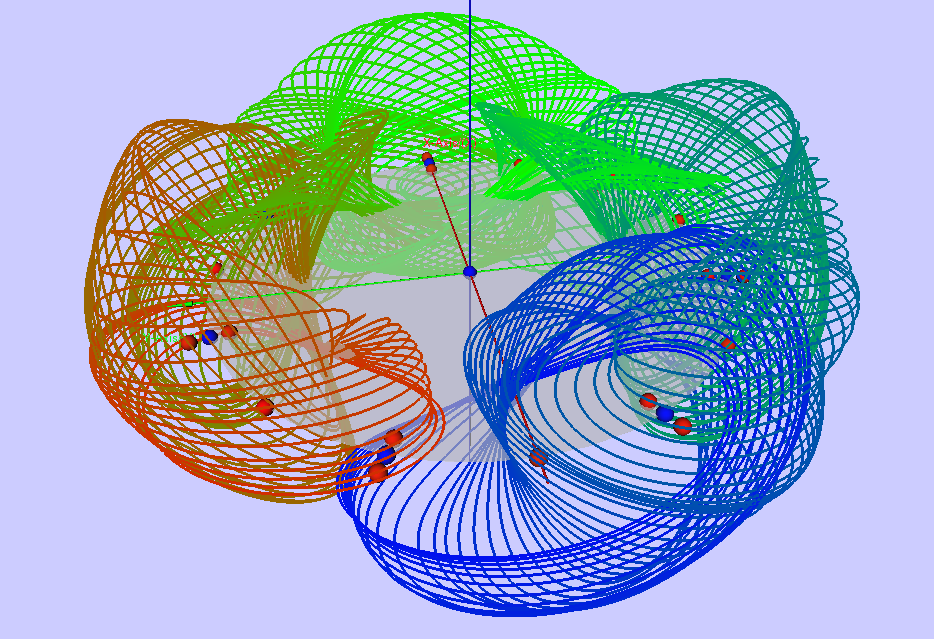}  \hskip1.9cm
\includegraphics{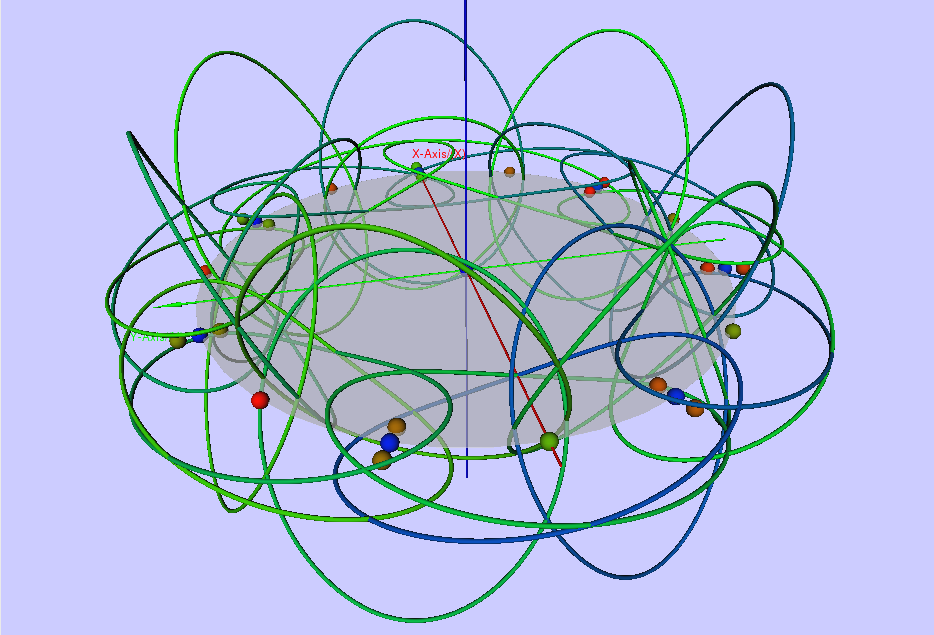} }
\end{center}
\caption{ Top Left: Part of the torus generated by the family $A_{3.1}$; the
torus closes up after making two full rounds. Top Right: Bifurcation orbits
along the torus generated by $A_{3.1}$. Bottom Left: Part of the torus
generated by the family $A_{3.2}$; This torus also closes up after making two
full rounds. Bottom Right: Bifurcation orbits along the torus generated by
$A_{3.2}$. }%
\label{fig:5}%
\end{figure}

\section{Secondary families}

\label{sec:second}

The bifurcations from the Lyapunov families coincide with the second breaking
of symmetries. The families that emanate from such bifurcation points (the red
lines in Figure~\ref{fig:2}) have three kinds of isotropy groups, each one
isomorphic to $\mathbb{Z}_{2}$. The number of symmetry-related branches is
equal to $2\times7$. Therefore, there are two such branches per equilibrium.

There is a symmetry-breaking from the planar families to solutions with
isotropy group
\begin{equation}
\mathbb{Z}_{2}(\kappa_{z})\text{.} \label{Se1}%
\end{equation}
These solutions have vertical component $z=0$. In this case the $\kappa_{y}%
$-symmetry is broken, so these planar orbits are asymmetric with respect to
the transformation that takes $y$ to $-y$. Such solutions are observed for the
third bifurcating family along $L_{2}$ and the fourth family that bifurcates
from $L_{1}$; see Figure~\ref{fig:4}.

The Halo families $H_{1}$ and $H_{2}$ bifurcate from the Lyapunov families
$L_{1}$ and $L_{2}$, respectively. Each of these has isotropy group
\begin{equation}
\mathbb{Z}_{2}(\kappa_{y})\text{,} \label{Se2}%
\end{equation}
\textit{i.e.}, their solutions have the property that $x(t)$ is even and
$y(t)$ is odd. Thus these spatial orbits are invariant under the
transformation that takes $y$ to $-y$.

The Axial families $A_{j}$, and similarly families $B_{j}$, bifurcate from
$V_{j}$, for $j=1,2,3$. Here the symmetry-breaking is from the group
$G_{V_{j}}$ to the group%

\begin{equation}
\mathbb{Z}_{2}(\pi\kappa_{y},\kappa_{z})\text{.} \label{Se3}%
\end{equation}
This means that the Axial solutions satisfy $\pmb{u}(t)=R_{z}R_{y}(-t+\pi)$
or, setting $\tilde{\pmb{u}}(t)=\pmb{u}(t+\pi/2)$, that
\[
\tilde{\pmb{u}}(t)=R_{z}R_{y}\tilde{\pmb{u}}(-t).
\]
Then $x(t)$ is even, and $y(t)$ and $z(t)$ are odd. Therefore these spatial
orbits are invariant under the transformation that takes $(y,z)$ to
$(-y,-z)$.

\section{Tertiary families}

\label{sec:third}

Tertiary symmetry-broken families (the green lines in Figure~\ref{fig:2})
correspond to families that bifurcate from solutions with isotropy groups
(\ref{Se1}), (\ref{Se2}), and (\ref{Se3}). Since these groups are isomorphic
to $\mathbb{Z}_{2}$, the tertiary symmetry-broken solutions have the trivial
isotropy group. Therefore, the number of symmetry-related branches of the
tertiary families is $4\times7$, i.e. four per equilibrium. In particular, the
symmetry-breaking bifurcations along the families $A_{3}$ and $B_{3}$ give
rise to families that generate invariant surfaces.

\begin{figure}[th]
\begin{center}
\resizebox{13cm}{!} {\includegraphics{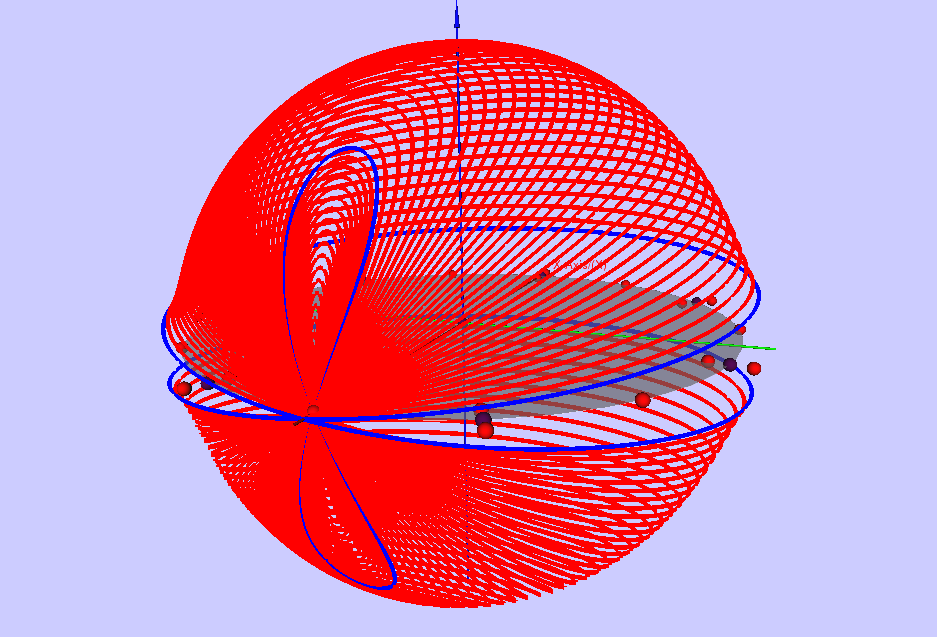}  \hskip1.9cm
\includegraphics{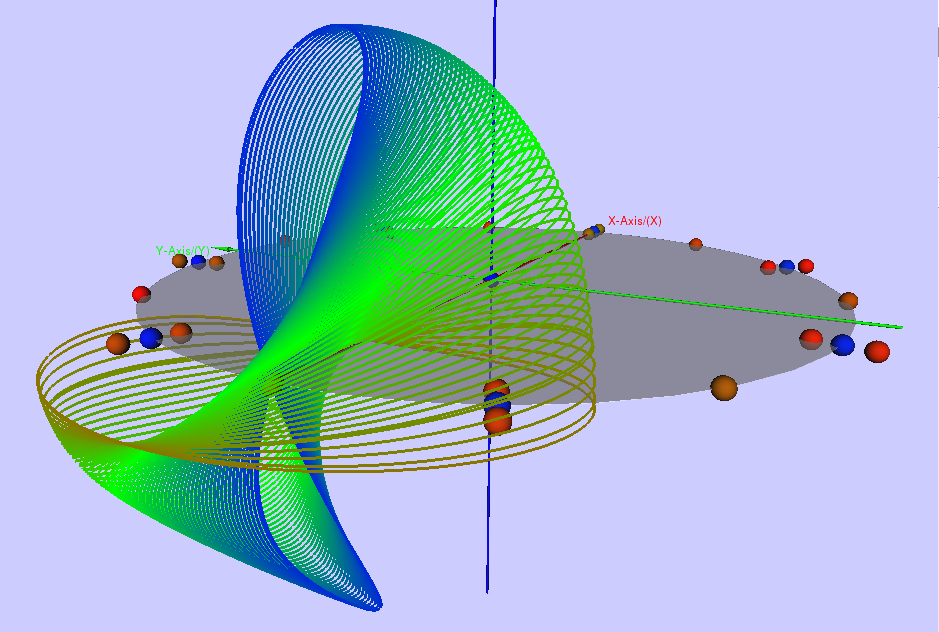} }
\par
\resizebox{13cm}{!} {\includegraphics{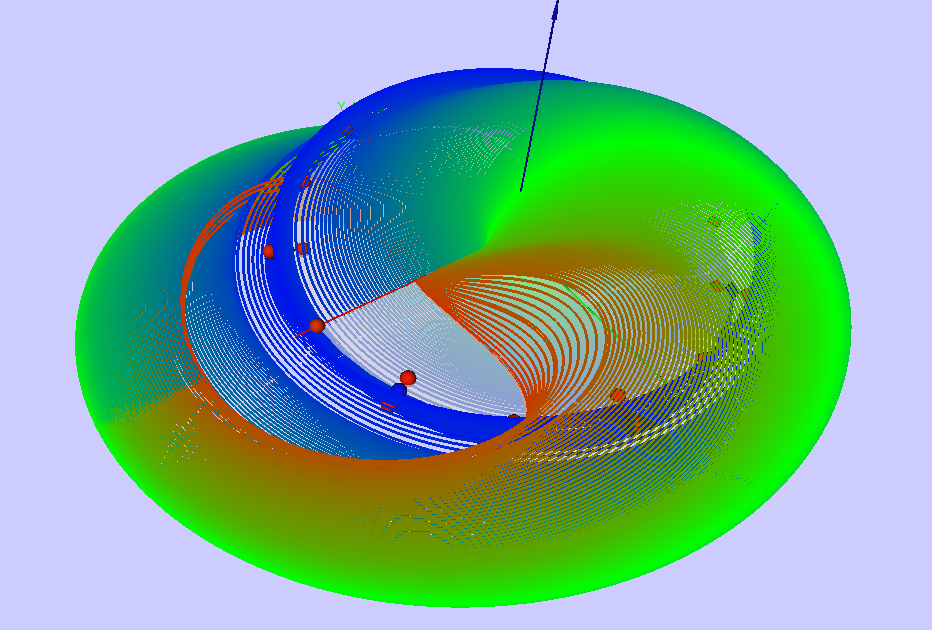}  \hskip1.9cm
\includegraphics{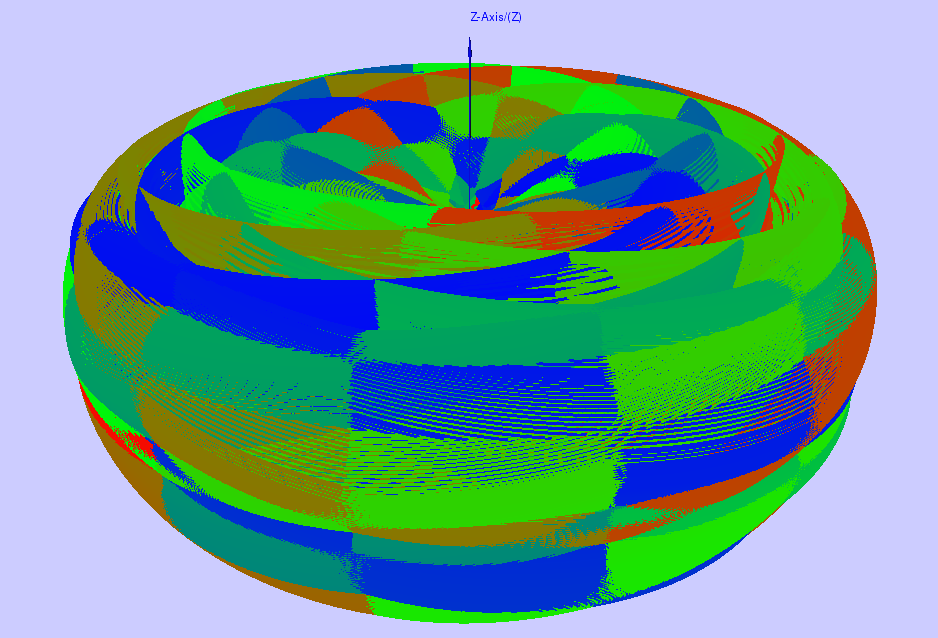} }
\end{center}
\caption{ Top Left: The Vertical family $V_{3}$, from the libration point
$\mathcal{L}_{3}$ until its second bifurcation orbit, with the bifurcation
orbits colored blue. Top Right: The Axial family $A_{3}$, from the orbit it
shares with $V_{3}$ and until the orbit it shares with the planar family
$P_{3}$. Bottom Left: The family $B_{3}$, from the orbit it shares with
$V_{3}$ and until the orbit it shares with the planar family $P_{3}$. Bottom
Right: One of several families that bifurcate from $B_{3}$. The orbits of this
family generate a torus that closes up after making two full rounds around the
$z$-axis. }%
\label{fig:6}%
\end{figure}

This is the case of the surface generated by the $A_{3.1}$ family in
Figure~\ref{fig:5} that reconnects to $A_{3}$ after a complete loop around the
central body, but to an orbit that is symmetric to the original one. Following
this surface, \textit{i.e.}, its orbits, for a second loop around the central
body, we obtain a double surface that interconnects the $7$ pairs of Axial
bifurcation orbits that emanate from each of the $7$ libration points
symmetry-related to $\mathcal{L}_{3}$. The surface generated by the orbits of
the family $A_{3.1}$ contains an additional $2\times7$ bifurcation orbits that
connect to Halo families $H_{2}$. Consequently there is a continuous path in
the bifurcation diagram between any of the $7$ symmetry-related libration
points $\mathcal{L}_{2}$ and $\mathcal{L}_{3}$ (see the bifurcation diagram
Figure~\ref{fig:2}).

The surface generated by the orbits of the family $A_{3.2}$ in
Figure~\ref{fig:5} is similar to that of $A_{3.1}$ that contains a set of
$2\times7$ bifurcation orbits along the Axial families $A_{3}$. These
bifurcation orbits are distinct from the bifurcation orbits along $A_{3}$ that
are interconnected via the family $A_{3.1}$. The family $A_{3.2}$ has
$2\times7$ extra bifurcation orbits that connect to Halo-like families that we
do not describe here.

Several families bifurcate from the family $B_{3}$ with trivial isotropy
group, one of which is illustrated in the bottom right panel of
Figure~\ref{fig:6}. This family is similar to $A_{3.1}$ and $A_{3.2}$ in that
it connects one Axial family, with isotropy group $\mathbb{Z}_{2}(\pi
\kappa_{y},\kappa_{z})$, to one Halo-like family, with isotropy group
$\mathbb{Z}_{2}(\kappa_{z})$. In future work we will report other families
bifurcating from $B_{3}$ that connect two Axial-like families with isotropy
group $\mathbb{Z}_{2}(\pi\kappa_{y},\kappa_{z})$.
\begin{figure}[th]
\begin{center}
\resizebox{13cm}{!} {\includegraphics{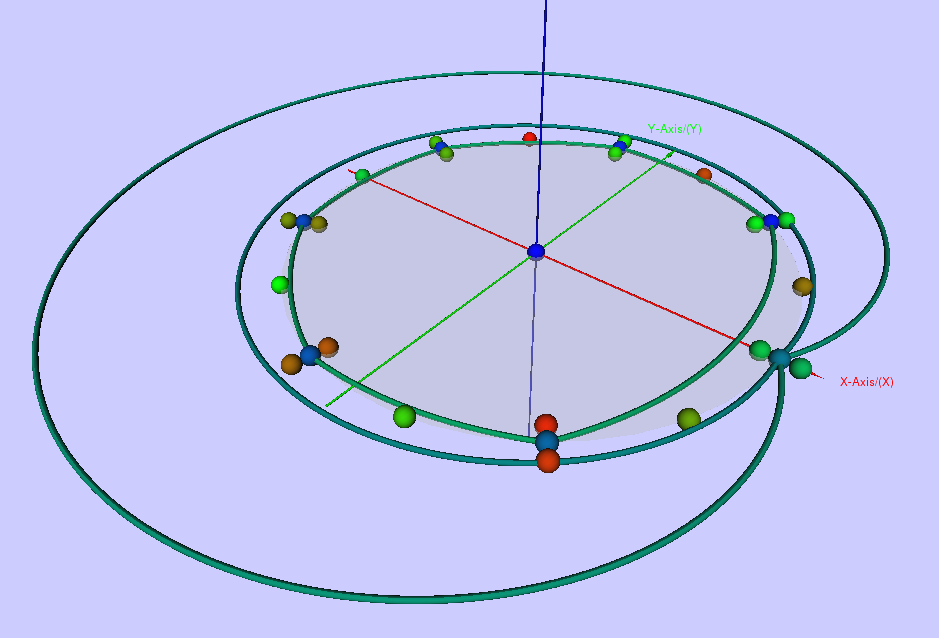}  \hskip1.9cm
\includegraphics{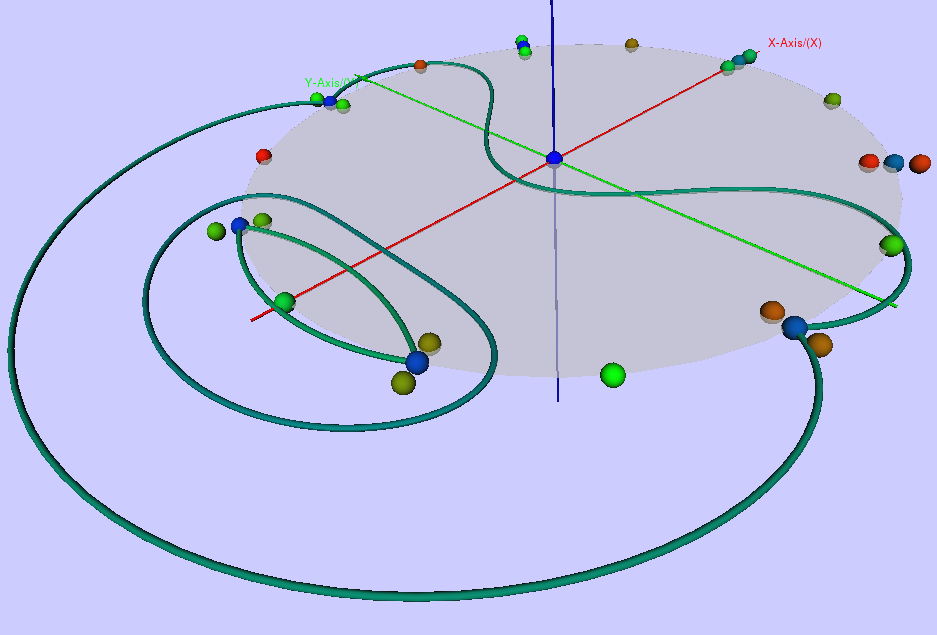} }
\end{center}
\caption{ Left: Three orbits from the planar, interplanetary family $C_{2}$
from which the vertical family $V_{2}$ bifurcates via a period-doubling.
Right: Three orbits from the planar, interplanetary family $P_{3}$ that
bifurcates from the Axial family $A_{3}$. }%
\label{fig:7}%
\end{figure}

\section{Interplanetary orbits}

\label{sec:four}

Several families connect to planar families that enclose more than one body
and that do not correspond to a planar Lyapunov family. Such families are
referred to as interplanetary in \cite{Ka08}, and indicated by black lines in
the bifurcation diagram of Figure~\ref{fig:2}.

The family $V_{2}$ has the symmetry group \eqref{GV} and connects to a planar
family $C_{2}$ in Figure~\ref{fig:7} via a reverse period-doubling
bifurcation, \textit{i.e.}, the vertical family $V_{2}$ arises from the planar
family via a period-doubling bifurcation. Indeed, the family $V_{2}$
bifurcates from the interplanetary family $C_{2}$ with isotropy group
\begin{equation}
\mathbb{Z}_{2}(\pi)\times\mathbb{Z}_{2}(\kappa_{y})\times\mathbb{Z}_{2}%
(\kappa_{z})\text{.} \label{Int}%
\end{equation}
These planar solutions are $\pi$-periodic and their orbits are invariant under
the transformation that takes $y$ to $-y$. This is consistent with the
symmetry-breaking phenomenon, since the isotropy group \eqref{GV} is contained
in the isotropy group of the interplanetary orbits (\ref{Int}).

The Lyapunov families $V_{1}$ and $V_{3}$ also connect to planar
interplanetary families, $C_{1}$ and $C_{3}$, respectively, via a
period-doubling bifurcation. Actually, $C_{1}$ and $C_{3}$ correspond to the
same family, \textit{i.e.}, $C_{1}=C_{3}$, and moreover, both $V_{1}$ and
$V_{3}$ bifurcate from exactly the same orbit along $C_{1}=C_{3}$, via a
degenerate period-doubling bifurcation with two Floquet multipliers at $-1$.

Another family that ends in an interplanetary family is $A_{3}$; see
Figure~\ref{fig:7}. Here the isotropy group of the planar family $P_{3}$ is
\[
\mathbb{Z}_{2}(\pi\kappa_{y})\times\mathbb{Z}_{2}(\kappa_{z})\text{.}%
\]
These solutions have the property that $x(t)$ is even, $y(t)$ is odd, and
$z=0$, so that the orbits are planar and symmetric with respect to the
transformation that takes $y$ to $-y$. Therefore the Axial family $A_{3}$,
with group $\mathbb{Z}_{2}(\pi\kappa_{y},\kappa_{z})$, corresponds to a
symmetry-breaking from the interplanetary family $P_{3}$.
\vskip0.25cm \textbf{Acknowledgements.} We thank Ramiro Chavez Tovar for his
assistance in producing the bifurcation diagram.



\begin{thebibliography}{99}                                                                                               %


\bibitem {BaEl04}D. Bang, B. Elmabsout. \emph{Restricted }$n+1$\emph{-body
problem: existence and stability of relative equilibria.} Celestial Mechanics
and Dynamical Astronomy, 89(4):305--318, 2004.

\bibitem {CaDo12}R. Calleja, E. Doedel, A. Humphries, A. Lemus-Rodr\'{\i}guez,
B. Oldeman. \emph{Boundary-value problem formulations for computing invariant
manifolds and connecting orbits in the circular restricted three-body
problem.} Celestial Mechanics and Dynamical Astronomy, 114(1-2):77-106, 2012.

\bibitem {GaIz10}C. Garc\'{\i}a-Azpeitia, J. Ize. \emph{Global bifurcation of
planar and spatial periodic solutions in the restricted }$\emph{n}$\emph{-body
problem.} Celestial Mechanics and Dynamical Astronomy, 110:217-227, 2011.

\bibitem {GaIz13}C.~Garc\'{\i}a-Azpeitia, J. Ize. \emph{Global bifurcation of
planar and spatial periodic solutions from the polygonal relative equilibria
for the } $n$\emph{-body problem.} J. Differential Equations,
254(5):2033--2075, 2013.

\bibitem {Henrard}J. Henrard. \emph{The web of periodic orbits at L4.}
Celestial Mech. Dynam. Astronom. 83(1):291--302, 2002.

\bibitem {Ka08}T.  Kalvouridis. \emph{Particle motions in Maxwell's ring
dynamical systems.} Celestial Mechanics and Dynamical Astronomy,
102(1-3):191--206, 2008.

\bibitem {LiMa}J. Llibre , R. Mart\'{\i}nez, C. Sim\'{o}. \emph{Tranversality
of the invariant manifolds associated to the Lyapunov family of periodic
orbits near L2 in the restricted three-body problem}. J. Differential
Equations, 58:104--156, 1985.

\bibitem {LoMa00}W. Koon, M. Lo, J. Marsden, S. Ross. \emph{Heteroclinic
connections between periodic orbits and resonance transitions in celestial
mechanics.} Chaos, 10(2):427-469, 2000.

\bibitem {Ma}J. Maxwell. On the stability of motions of Saturns rings.
Macmillan and Co., Cambridge, 1859.

\bibitem {Moekel92}R. Moeckel. \emph{Linear stability of relative equilibria
with a dominant mass}. J. of Dynamics and Differential Equations, 6:37--51, 1994.

\bibitem {VaKo07}R. Vanderbei, E. Kolemen. \emph{Linear stability of ring
systems. }The Astronomical Journal, 133:656--664, 2007.

\bibitem {}A. Vanderbauwhede. \emph{Branching of Periodic Orbits in
Hamiltonian and Reversible Systems. }Equadiff 9: Proceedings of the 9th
conference, Brno, 169-181, 1997.

\bibitem {Ro00}G. Roberts. \emph{Linear stability in the } $1+n$\emph{-gon
relative equilibrium.} In J.~Delgado, editor, Hamiltonian Systems and
Celestial Mechanics. HAMSYS-98. World Sci. Monogr. Ser. Math. 6, 303--330.
World Scientific, 2000.

\bibitem {DoVa}E. Doedel, R. Paffenroth, H. Keller, D. Dichmann, J.
Gal\'{a}n-Vioque, A. Vanderbauwhede. \emph{Continuation of periodic solutions
in conservative systems with application to the }$\emph{3}$\emph{-body
problem.} Int. J. Bifurcation Chaos Appl. Sci. Eng. 13:1--29, 2003.
\end{thebibliography}
\end{document}